\documentclass[12pt]{amsart}
\usepackage{amsmath, amsfonts, latexsym, array, amssymb, }
\usepackage{psfrag}
\usepackage[all]{xypic}
\usepackage{graphicx}
\newtheorem{thm}{Theorem}[section]

\newtheorem{lem}[thm]{Lemma}
\newtheorem{prop}[thm]{Proposition}

\newtheorem{que}[thm]{Question}
\newtheorem{defn}[thm]{Definition}

\def\eps{\varepsilon}

\newcommand{\occult}[1]{}

\newcommand{\new}[1]{{\bf #1}}

\begin{document}

\title[Intrinsic ergodicity] {Intrinsic ergodicity for certain nonhyperbolic robustly transitive systems}
\author{J\'er\^ome Buzzi, Todd Fisher}
\address{C.N.R.S. \& D\'epartement de Math\'ematiques, Universit\'e Paris-Sud, 91405 Orsay, France}
\address{Department of Mathematics, Brigham Young University, Provo, UT 84602}
\email{jerome.buzzi@math.u-psud.fr}
\email{tfisher@math.byu.edu}
\thanks{}

\date{June, 19, 2008}
\commby{}

\begin{abstract}
We show that a class of robustly transitive diffeomorphisms originally described by Ma\~{n}\'{e} are intrinsically ergodic. More precisely, we obtain an open set of diffeomorphisms which fail to be uniformly hyperbolic and structurally stable, but nevertheless have constant entropy and isomorphic unique measures of maximal entropy.
 \end{abstract}

\maketitle

\section{Introduction}

Let $f$ be a diffeomorphism of a manifold $M$ to itself.  The diffeomorphism $f$ is {\it transitive} if there exists a point $x\in M$ where
$$\mathcal{O}^+_f(x)=\{f^n(x)|n\in\mathbb{N}\}$$
is dense in $M$.
It is {\it robustly transitive} \cite[Ch. 7]{BDV} if there exists a neighborhood $U$ of $f$ in the space $\mathrm{Diff}^1(M)$ of $C^1$ diffeomorphisms  such that each $g$ in $U$ is transitive. 
Since robust transitivity is an open condition, it is an important component of the global picture of dynamical systems~\cite{Pal}.   

The first examples of robustly transitive diffeomorphisms were transitive Anosov diffeomorphisms:  recall that a diffeomorphism is Anosov if the entire manifold is a hyperbolic set under the action of the diffeomorphism.  Nonhyperbolic robustly transitive diffeomorphisms were first constructed by Shub~\cite{Shu1} and Ma\~{n}\'{e}~\cite{man78}.
These examples satisfy a weaker hyperbolic condition called partial hyperbolicity 
(see Sec. \ref{sec:background}).
It is interesting to note when results for Anosov diffeomorphisms continue to hold and when the properties are very different.  For instance, $C^1$-structural stability holds for Axiom A systems with strong transversality and no others~\cite{man78}.  In this paper we analyze measures of maximal entropy and a related notion of stability for some class of non-Anosov robustly transitive diffeomorphisms based on Ma\~{n}\'e's example.

To state our results, we need to give some definitions.  Dynamical entropies are measures of the complexity of orbit structures \cite{Bow}. The topological entropy, $h_{\mathrm{top}}(f)$, considers all the orbits, whereas the measure theoretic entropy, $h_{\mu}(f)$, focuses on those ``relevant" to a given invariant probability measure $\mu$.  The variational principle, see for example~\cite[p. 181]{KH1}, says that if $f$ is a continuous self-map of a compact metrizable space and $\mathcal{M}(f)$ is the set of invariant probability measures for $f$, then
$$h_{\mathrm{top}}(f)=\sup_{\mu\in\mathcal{M}(f)}h_{\mu}(f).$$  A measure $\mu\in\mathcal{M}(f)$ such that $h_{\mathrm{top}}(f)=h_{\mu}(f)$ is a {\it measure of maximal entropy}. By a theorem of Newhouse \cite{Newhouse} $C^\infty$ smoothness implies the existence of such measures (but finite smoothness does not according to Misiurewicz \cite{Misiurewicz}). If there is a unique measure of maximal entropy, then $f$ is called {\it intrinsically ergodic}.

\begin{defn}
We say $f\in\mathrm{Diff}^1(M)$ is \new{intrinsically stable} if there exists a neighborhood $U$ of $f$ such that each $g$ in $U$ has a unique measure of maximal entropy $\mu_g$ and
all $\mu_g$ define isomorphic measure-preserving transformations.
\end{defn}

Newhouse and Young~\cite{NY} have shown that the robustly transitive diffeomorphisms constructed by Shub on $\mathbb T^4$ are intrinsically stable (and in particular intrinsically ergodic).   The present work extends this to the robustly transitive diffeomorphisms constructed by Ma\~{n}\'{e} on $\mathbb T^3$.

\begin{thm}\label{t.mane}
For any $d\geq 3$, there exists a non-empty open set $U$ 
in $\mathrm{Diff}(\mathbb T^d)$ satisfying:
 \begin{itemize}
  \item each $f\in U$ is strongly partially hyperbolic, robustly transitive, and intrinsically stable 
(in particular the topological entropy is locally constant at $f$);
  \item no $f\in U$ is Anosov or structurally stable.
 \end{itemize}
\end{thm}

This raises the following question.

\begin{que} Is every robustly transitive diffeomorphism intrinsically ergodic? intrinsically stable?
\end{que}

We note that examples of Kan \cite{BDV, Kan} suggest that the answer might be negative. These are robustly transitive systems within $C^1$ self-maps of the compact cylinder preserving the boundary which admits two SRB measures on the boundary that are also measures of maximal entropy.

In dimension three we know that every robustly transitive system is partially hyperbolic~\cite[p. 128]{BDV}.  The added structure of partial hyperbolicity could help solve the above question in the affirmative for 3-manifolds.

In a follow up paper~\cite{BF1} we analyze a set of robustly transitive diffeomorphisms on $\mathbb{T}^4$, based on examples of Bonatti and Viana~\cite{BV} which have the weakest possible form of hyperbolicity for robustly transitive diffeomorphisms: dominated splitting~\cite{BDP, DPU} (see Sec. \ref{sec:background}).

We note that Hua, Saghin, and Xia \cite{Hua} have also proved local constancy of the topological entropy, for instance in the case of partially hyperbolic diffeomorphisms $C^1$ close to toral automorphisms with at most one eigenvalue on the unit circle have locally constant entropy.
 
Just before submitting this paper to the arxiv we received a communication from Carlos Matheus about a work in progress on the existence of maximal measures and more generally equilibrium states for such systems.  Afterwards, we received a manuscript from Sambarino and Vasquez~\cite{SV} similar to the present work.

\medbreak
\noindent{\bf Acknowledgment.} We thank Jean-Paul Allouche and Lennard Bakker for their help on Pisot numbers, and Sheldon Newhouse for helpful discussions.

\section{Background}\label{sec:background}

We now review a few facts on entropy, hyperbolicity, and partial hyperbolicity.

Let $X$ be a compact metric space and $f$ be a continuous self-map of $X$.  Fix $\epsilon>0$ and $n\in\mathbb{N}$.  Let $\mathrm{cov}(n,\epsilon, f)$ be the minimum cardinality of a covering of $X$ by $(\epsilon,n)$-balls, i.e., sets of the form 
$$\{y\in X:\; d(f^k(y),f^k(x))<\epsilon \textrm{ for all }0\leq k\leq n
\}.$$
The {\it topological entropy} is \cite{Bow}
$$h_{\mathrm{top}}(f)=\lim_{\epsilon\rightarrow 0}(\limsup_{n\rightarrow\infty}\frac{1}{n}\log \mathrm{cov}(n,\epsilon, f)).$$

Let $Y\subset X$ and $\mathrm{cov}(n,\epsilon,f,Y)$ be the minimum cardinality of a cover of $Y$ by $(n,\epsilon)$-balls.  Then the {\it topological entropy of $Y$ with respect to $f$} is
$$
{h}_{\mathrm{top}}(f,Y)=\lim_{\epsilon\rightarrow 0}\limsup_{n\rightarrow\infty}\frac{1}{n}\log \mathrm{cov}(n,\epsilon,f,Y).$$
If $(X,f)$ and $(Y,g)$ are continuous and compact systems and $\phi:X\rightarrow Y$ is a continuous surjection such that $\phi\circ f=g\circ \phi$, then $h_{\mathrm{top}}(g)\leq h_{\mathrm{top}}(f)$ ($f$ is called an {\it extension} of $g$ and $g$ is called a {\it factor} of $f$).  For the definition of measure theoretic entropy refer to~\cite[p. 169]{KH1}.

An invariant set $\Lambda$ is {\it hyperbolic} for $f\in\mathrm{Diff}(M)$ if  there exists an invariant splitting $T_{\Lambda}M=E^s\oplus E^u$ and an integer $n\geq1$ such that $Df^n$ uniformly contracts $E^s$ and uniformly expands $E^u$: for any point
$x\in\Lambda$,
$$
\begin{array}{llll}
\| Df_{x}^{n}v\|\leq \frac12\| v\|,\textrm{ for }v\in
E^{s}_x \textrm{, and}\\
\| Df_{x}^{-n}v\|\leq \frac12\| v\|,\textrm{ for }v\in
E^{u}_x. \end{array}
$$

If $A\in \mathrm{GL}(d,\mathbb{Z})$ has no eigenvalues on the unit circle, then the induced map $f_A$ of the $d$-torus is called a {\it hyperbolic toral automorphism}.  By construction any hyperbolic toral automorphism is Anosov.

If $\Lambda$ is a hyperbolic set, $x\in\Lambda$, and $\epsilon>0$ sufficiently small, then  the
\textit{local stable and unstable manifolds} at $x$ are respectively:
$$
\begin{array}{llll}
W_{\epsilon}^{s}(x,f)=\{ y\in M\, |\textrm{ for all }
n\in\mathbb{N}, d(f^{n}(x), f^{n}(y))\leq\epsilon\},\textrm{
and}\\
W_{\epsilon}^{u}(x,f)=\{ y\in M\, |\textrm{ for all }
n\in\mathbb{N}, d(f^{-n}(x), f^{-n}(y))\leq\epsilon\}.
\end{array}$$
The \textit{stable and unstable manifolds} of $x$ are respectively:
\begin{align*}
   W^s(x,f)&=\{y\in M\, | \lim_{n\to\infty} d(f^n(y),f^n(x))=0 \}, \textrm{ and}\\
   W^u(x,f)&=\{y\in M\, | \lim_{n\to\infty} d(f^{-n}(y),f^{-n}(x))=0\}.
\end{align*}
They can be obtained from the local manifolds as follows:
 $$
 \begin{array}{llll}
W^s(x,f)=\bigcup_{n\geq 0}f^{-n}\left(
W_{\epsilon}^s(f^n(x),f)\right),\textrm{ and}\\
W^u(x,f)=\bigcup_{n\geq
0}f^{n}\left(W_{\epsilon}^u(f^{-n}(x),f)\right).
\end{array}$$
For a $C^r$ diffeomorphism the stable and unstable manifolds of a hyperbolic set are $C^r$
injectively immersed submanifolds.

An {\it $\epsilon$-chain} from a point $x$ to a
point $y$ for a diffeomorphism $f$ is a sequence $\{x=x_0,...,x_n=y\}$ such
that 
$$d(f(x_{j-1}),x_j)<\epsilon\textrm{ for all }1\leq j\leq n.$$  A standard result that applies to Anosov diffeomorphisms is the
Shadowing Theorem, see for example~\cite[p. 415]{Rob}.  Let $\{x_j\}_{j=j_1}^{j_2}$
be an $\epsilon$-chain for $f$.  A point $y$ {\it
$\delta$-shadows} $\{x_j\}_{j=j_1}^{j_2}$ provided
$d(f^j(y),x_j)<\delta$ for $j_1\leq j\leq j_2$.  We remark that there are much more general versions of the next theorem, but the following statement will be sufficient for the present work.

\begin{thm}(Shadowing Theorem) If $f$ is an Anosov diffeomorphism, then given any $\delta>0$ sufficiently small there exists an $\epsilon>0$  such
that if $\{ x_j\}_{j=j_1}^{j_2}$ is an $\epsilon$-chain for $f$, then there is a $y$ which $\delta$-shadows
$\{x_j\}_{j=j_1}^{j_2}$.  If $j_2=-j_1=\infty$, then $y$ is unique. If, moreover, the $\epsilon$-chain is periodic, then $y$ is periodic.
\end{thm}

A diffeomorphism $f:M\rightarrow M$ has a {\it dominated splitting} if there exists an invariant  splitting $TM=E_1\oplus\cdots E_k$, $k\geq 2$, (with no trivial subbundle) and  an integer $l\geq1$ such that for each $x\in M$, $i<j$, and unit vectors $u\in E_i(x)$ and $v\in E_j(x)$, one has
$$\frac{\|Df^l(x)u\|}{\|Df^l(x)v\|}<\frac{1}{2}.$$
A diffeomorphism $f$ is {\it partially hyperbolic} if there is a dominated splitting $TM=E_1\oplus\cdots\oplus E_k$ and $n\geq 1$ such that $Df^n$ either uniformly contracts $E_1$ or uniformly expands $E_k$.  We say $f$ is {\it strongly partially hyperbolic} if there exists a dominated splitting $TM=E^s\oplus E^c\oplus E^u$ and $n\geq 1$ such that $Df^n$ uniformly contracts $E^s$ and uniformly expands $E^u$.

For $f$ a strongly partially hyperbolic diffeomorphism we know there exist unique families $\mathcal{F}^u$ and $\mathcal{F}^s$ of injectively immersed submanifolds such that $\mathcal{F}^i(x)$ is tangent to $E^i$ for $i=s,u$, and the families are invariant under $f$, see~\cite{HPS}.  These are called, respectively,  the unstable and stable laminations\footnote{A $C^f$ foliation is a partition of the manifolds locally $C^r$-diffeomorphic (or homeomorphic if $f=0$) to a partition of $\mathbb R^d$ into $k$-planes for some $0\leq k\leq d$. A lamination is a $C^0$ foliation with $C^1$ leaves.}
 of $f$.  For the center direction, however,
there are examples where there is no center lamination~\cite{Wil}.  For a strongly partially hyperbolic diffeomorphism with a 1-dimensional center bundle
it is not known if there is always a lamination tangent to the center bundle, and that if there is a $C^1$ center foliation, then it is structurally stable~\cite{HPS}. Let us quote a special case of this result:

\begin{thm}\label{t.HPS}\cite[Theorems (7.1) and (7.2)]{HPS}
Let $f$ be a $C^1$ diffeomorphism of a compact manifold $M$. If $f$ is strongly partially hyperbolic
with a $C^1$ central foliation $\mathcal F$, then any $g$ $C^1$-close to $f$ also has a $C^1$ central lamination $\mathcal G$ and there is a homeomorphism $h:M\to M$ such that for all $x\in M$, (i) the leaf $\mathcal F_x$ is mapped by $h$ to the leaf $\mathcal G_{hx}$; (ii) $g(\mathcal G_{hx})= \mathcal G_{h(fx)}$.
\end{thm}

This applies in particular to the Ma\~{n}\'{e} example.

\section{Intrinsic ergodicity for Ma\~{n}\'{e}'s robustly transitive diffeomorphisms}

Ma\~{n}\'{e}'s example of a robustly transitive dynamical system that is not Anosov was constructed on $\mathbb{T}^3$.  We will use his construction for diffeomorphisms of higher dimensional tori.

\begin{figure}[htb]
\begin{center}
\psfrag{A}{$f_A$}
\psfrag{B}{$f_0$}
\psfrag{q}{$q$}
\psfrag{r}{$q_1$}
\psfrag{s}{$q$}
\psfrag{t}{$q_2$}
\includegraphics{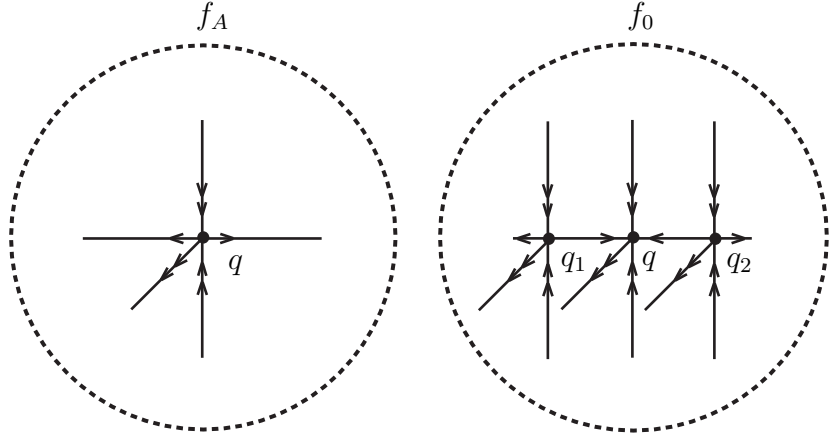}
\caption{Ma\~{n}\'{e}'s construction}\label{f.mane}
\end{center}
\end{figure}

We fix some dimension $d\geq3$ and let $A\in\mathrm{GL}(d,\mathbb{Z})$ be a hyperbolic toral automorphism with only one eigenvalue inside the unit circle and all eigenvalues real, positive, simple, and irrational.
Let $\lambda_s$ be the unique modulus less than $1$ and $\lambda_c$ be the smallest
of the moduli greater than $1$. 

We denote the induced linear Anosov system on $\mathbb{T}^d$ by $f_A$ and let
$\mathcal{F}^c$ be the foliation corresponding to the eigenvalue $\lambda_c$; so locally at each point $\mathcal{F}^c$ is just a line segment in the direction of the eigenvector associated with $\lambda_c$.
Similarly, $\mathcal{F}^s$ and $\mathcal F^u$ are the foliations corresponding to the eigenvalue $\lambda_s$ and all the eigenvalues greater than $\lambda_c$, respectively.   Since all eigenvalues are irrational, each leaf of $\mathcal{F}^s$, $\mathcal F^c$, and $\mathcal{F}^u$ is dense in $\mathbb{T}^d$.

Such matrices can be built for any $d\geq3$ as companion matrices to the minimal polynomial over $\mathbb Q$ of a Pisot number whose algebraic conjugates are all real. Such numbers are given by Theorem 5.2.2 in~\cite[p. 85]{Bertin} (the proof implies that the conjugates are real). The moduli are then pairwise distinct by \cite{Smyth}.

Without loss of generality, we may assume that $f_A$ has at least two fixed points and that any unstable eigenvalue other than $\lambda_c$ has modulus greater than $3$ (if not, replace $A$ by some power).

Let $p$ and $q$ be fixed points under the action of $f_A$ and $\rho>0$ be a small number to be determined below. Following the construction in~\cite{man78} we define $f_0$ by modifying $f_A$ in a sufficiently small domain $C$ contained in $B_{\rho/2}(q)$ keeping invariant the foliation $\mathcal{F}^c$.  So there is a neighborhood $U$ of $p$ such that $f_A|_U=f_0|_U$.  Inside $C$ the fixed point $q$ undergoes a pitchfork bifurcation in the direction of the foliation $\mathcal{F}^c$.  The stable index of $q$ increases by $1$, and two other saddle points with the same stable index as the initial $q$ are created. (See Figure~\ref{f.mane}.)

The resulting diffeomorphism $f_0$ is strongly partially hyperbolic with a $C^1$ center foliation $\mathcal{F}^c$. 
According to~\cite{man78}, it is also robustly transitive (in fact topologically mixing~\cite[p. 184]{BDV}) for $\rho>0$ sufficiently small.

The next proposition will be helpful in the proof of Theorem~\ref{t.mane}.

\begin{prop}\label{p.semi}(Shadowing proposition)
Let $f_A$ be an Anosov diffeomorphism of the $d$-torus, $d\geq3$, as above.
Let $f\in\mathrm{Diff}^1(\mathbb{T}^d)$ satisfy the following properties:
\begin{enumerate}
\item[(a)] $f$ contains a fixed point $p \in\mathbb{T}^d$ with $\overline{W^s(p)}=M$,
\item[(b)] there exist constants $\epsilon>0$ and $\delta>0$ such that each $\epsilon$-chain under $f_A$ is $\delta$-shadowed by an orbit under $f_A$ and $3\delta$ is an expansive constant for $f_A$, (i.~e. if $x,y\in\mathbb{T}^d$ and  $d(f^n_A(x), f^n_A(y))<3\delta$ for all $n\in\mathbb{Z}$, then $x=y$),  and
\item[(c)] each $f$-orbit is an $\epsilon$-chain for $f_A$.
\end{enumerate}
Then the map $\pi:\mathbb{T}^d\rightarrow \mathbb{T}^d$, where $\pi(x)$ is the point in $\mathbb{T}^d$ that under the action of $f_A$ will $\delta$-shadow the $f$-orbit of $x$, is a semiconjugacy from $f$ to $f_A$, i.e., it is a continuous and onto map with $\pi\circ f=f_A\circ\pi$.
\end{prop}
\noindent{\bf Proof.} By the shadowing theorem we know that the map $\pi$ is well-defined and that $\pi (f(x))=f_A(\pi (x))$ and $d(\pi(x),x)<\delta$. We need to see that $\pi$ is continuous~\cite[Theorem 7.8]{Shu1} and surjective. It is probably folklore, but we provide a proof for the convenience of the reader.  

To show that $\pi$ is continuous  we take a sequence $x_n\rightarrow x$ and show that $\pi (x_n)\rightarrow \pi(x)$.  Fix $M\in\mathbb{N}$.  Then there exists an $N(M)\in\mathbb{N}$ such that for each $n\geq N(M)$
$$d(f^j(x_n), f^j(x))<\delta\textrm{ for all }-M\leq j\leq M.$$
We then have
$$d(f_A^j(\pi(x_n)), f_A^j(\pi(x)))<3\delta\textrm{ for all }-M\leq j\leq M$$ where $n\geq N(M)$. It follows that for any limit point $y$ of the sequence $\{\pi(x_n)\}$ we have
 \begin{equation}\label{eq:shadlim}
   d(f_A^j(y), f_A^j(\pi (x)))\leq 3\delta\textrm{ for all }j\in\mathbb{Z}.
 \end{equation}
Since $3\delta$ is an expansive constant for $f_A$ this implies that $y=\pi (x)$ and $\pi (x_n)$ converges to $\pi (x)$.

We now show that $\pi$ is surjective.  Let $x\in W^s(p)$ with $d(x,p)>2\delta$. $f^n(x)\to p$ and $\pi(p)=p_A$, hence $\pi(x)\in W^s(p_A)$. Also 
$$d(\pi(x),p_A)>d(x,p)-2\delta>0.$$ 
Thus the segment $[\pi(x),\pi(f(x)))_s$ along $W^s(p_A)$ is non-trivial. By continuity of $\pi$ we know that 
 $$\pi([x,f(x))_s)\supset[\pi(x),f(\pi(x)))_s.$$ 
 It follows that the image of $\pi$ contains one of the connected components of $W^s(p_A)\setminus\{p_A\}$. Hence the image of $\pi$ is dense.  As $\pi$ is continuous we know that $\pi$ is surjective. $\Box$

\medbreak

We remark that the surjectivity of the map $\pi$ in Proposition~\ref{p.semi} can also be obtained by a topological argument.  Namely, that $\pi$ depends continuously on $f$ and is the identity for $f=f_A$.  On the torus this forces the surjectivity for $f$ homotopic to $f_A$.

We shall also use the following (folklore) fact:
\begin{lem}
Let $g:\mathbb T^d\to \mathbb T^d$ be an injective continuous self-map.
Let $K$ be a compact curve such that the lengths of all its iterates, $g^n(K)$, $n\geq0$,
are bounded by a constant $L$. Then $h(g,K)=0$.
\end{lem}

\noindent{\bf Proof of Lemma}
For each $n\geq0$, there exists a subset $K(\eps,n)$ of $g^n(K)$ with cardinality at most $L/\eps+1$ dividing $g^n(K)$ into curves with length at most $\eps$. Observe that $\bigcup_{0\leq k<n} g^{-k}K(\eps,k)$ is an $(n, \epsilon)$-cover of $K$ with subexponential cardinality. $\Box$

\medbreak

\noindent{\bf Proof of Theorem~\ref{t.mane}} 
The strategy of the proof of Theorem~\ref{t.mane} is to use the semiconjugacy $\pi_g$ from Proposition~\ref{p.semi} and to show that for each $x\in\mathbb{T}^d$ and each $g$ $C^1$-close to $f_0$, the set $\pi_g^{-1}(x)$ is a compact interval of bounded length contained in a center leaf, and $\pi^{-1}_g(x)$ is a unique point for almost every $x$.
We note that the measure of maximal entropy for $f_A$ is Lebesgue measure, denoted $\mu$, on $\mathbb{T}^d$.

We claim that for $\rho>0$ small enough, any diffeomorphism $f$ that is $C^1$ close to the previously constructed diffeomorphism $f_0$, satisfies the hypothesis of Proposition~\ref{p.semi}. Hypothesis (b) and (c) are clear. Let us show (a).

By Theorem \ref{t.HPS}, there is a neighborhood $U_0$ of $f_0$ such that each $g\in U_0$ is strongly partially hyperbolic with a center lamination $\mathcal{F}^c_g$ close to that of the center foliation $\mathcal{F}^c$. In particular they both have dimension $1$ with bounded ``curvature", for any $g\in U$: if $x,y,z$ are on the same central leaf in that order with $x,y\in B(z,2\delta)$ then $d(z,y)<d(z,x)$.

To show (a) we let $V\ne\emptyset$ be an open set in $\mathbb{T}^d$ and let $\sigma\subset V$ be a connected piece of a center leaf. By density of the whole leaf, $\sigma$ is eventually expanded to become $\delta$-dense for any $\delta>0$.
Let 
$$D^{su}_{\epsilon}(p)=\bigcup_{y\in W^s_{\epsilon}(p)}W^{uu}_{\epsilon}(p)$$ 
where $W^{uu}_{\epsilon}(p)$ is the connected component of $\mathcal{F}^u(y)\cap B_{\epsilon}(p)$ containing $y$.
The set $D^{su}_{\epsilon}(p)$ is transverse to the center direction.  Therefore, there exists an arbitrarily large $n\geq 0$ such that $f^n(\sigma)\cap D^{su}_{\epsilon}(p)\neq\emptyset$.  Hence, there exists some $y\in W^s_{\epsilon}(p)$ with $W^{uu}_{\epsilon}(y)\cap f^n(\sigma)\neq\emptyset$.  As $\mathcal{F}^u$ is uniformly contracted under $f^{-1}$ this implies that $f^{-n}(y)\in W^s(p)\cap V$ and (a) follows.

 Let $r>0$ be an expansive constant for $f_A$ and fix a neighborhood $U\subset U_0$ of $f_0$ such that each $g\in U$ satisfies the hypothesis of Proposition~\ref{p.semi} with $0<\epsilon<\delta<\min(r/3,\rho)$.  For each $g\in U$ we denote $\pi_g$ as the semiconjugacy mapping $g$ to $f_A$ given by Proposition~\ref{p.semi}.

Let $\mu$ be Lebesgue measure on $\mathbb{T}^d$ and set
\begin{equation}\label{eq.m}
m=\mu(B(q, 3\rho))>0.
\end{equation}
The above construction implies is such that the maximum contraction in the center direction, denoted $b(f)$, satisfies
\begin{equation}\label{eq.2}
\lambda_c^{1-m}b(f)^{2m}>1
\end{equation}
where $m$ is defined in (\ref{eq.m}).

Fix $\gamma>0$ such that $(\lambda_c-\gamma)^{1-m}(b(f)-\gamma)^{2m}>1$
. Possibly by reducing $U$, we may and do assume that 
$d_{C^1}(f_0,g)<\gamma$ and that robust transitivity holds for all $g\in U$.

Fix $g\in U$ and suppose that $y_1,y_2\in \pi_g^{-1}(x)$. By construction of $\pi_g$,
this implies $d(g^n(y_1),g^n(y_2))<2\delta$ for all $n\in\mathbb Z$. The normal hyperbolicity of the center lamination implies that such $y_1$ and $y_2$ must lie in the same center leaf.
By the bounded curvature property, the whole segment of $\mathcal F^c$ between $y_1$ and $y_2$
stays within $2\delta<r$ of the orbit of $y_1$, hence its image by $\pi_g$ stays within $\epsilon +2\delta<r$ of the orbit of $x$ so this interval must be contained in $\pi_g^{-1}(x)$.  It follows that the set $\pi_g^{-1}(x)$ is a compact interval in a center leaf which keeps
a bounded length under all iterates of $g$. The above lemma implies that $h(g,\pi_g^{-1}(x))=0$ for all $x\in\mathbb{T}^d$.

We now show that the topological entropy is constant in $U$.  For $g\in U$ we know that $f_A$ is a topological factor of $g$.  This implies that $h(f_A)\leq h(g)$.
In~\cite{Bow} Bowen shows that
  $$h(g)\leq h(f_A) +\sup_{x\in\mathbb{T}^d}h(g,\pi_g^{-1}x).$$
The last entropy is zero, hence the diffeomorphisms $f_A$ and $g$ have equal topological entropy.

Let $\mathcal{M}(g)$ be the collection of Borel invariant probability measures for $g$.
From the Hahn-Banach theorem we know that there exists an invariant measure $\bar{\mu}$ such that
$(\pi_g)_*\bar{\mu}=\mu$.  Since $g$ is an extension of $f_A$ we know that
$h_{\bar\mu}(g)\geq h_{\mu}(f_A)= h(f_A)=h(g)$: $\bar{\mu}$ is a measure of maximal entropy for $g$.

Now take $\nu$ an arbitrary measure of maximal entropy for $g$ and let us show that $\nu=\bar\mu$.  From results of Ledrappier and Walters  in~\cite{LW} we know that
$$
h_{\nu}(g) = h_{(\pi_g)_*\nu}(f_A) + \int_{\mathbb{T}^d} h(g,\pi_g^{-1}x)d(\pi_g)_*\nu
 = h_{(\pi_g)_*\nu}(f_A).
$$
The intrinsic ergodicity of $f_A$ implies that $(\pi_g)_*\nu=\mu$.

To prove that $g$ itself is intrinsically ergodic we show that $\pi_g$ is almost everywhere one-to-one, i.e. that Lebesgue almost every point in $\mathbb{T}^d$ has a unique pre-image under $\pi_g$.  Since $\mu$ is ergodic for $f_A$ we know from Birkhoff's ergodic theorem (see~\cite[p. 274]{Rob}) that for $\mu$-almost every $x\in\mathbb{T}^d$ we have
\begin{equation}\label{eq.3}
\lim_{n\rightarrow \infty}\frac{1}{n}\sum_{i=1}^n\chi_{B(q, \rho+2\delta)}(f_A^i(x))=\mu(B(q,\rho+2\delta))=m.
\end{equation}

Fix $g\in U$ and let
$$a(g)=\min_{x\in \mathbb{T}^3-B(q, \rho)}Dg_x \mathcal{F}^c(x)\geq \lambda_c-\gamma$$
and
$$b(g)=\min_{x\in B(q, \rho)} Dg_x\mathcal{F}^c(x)\geq b(f)-\gamma.$$
So $a(g)$ measures the minimum expansion in $\mathbb{T}^d-B(q,\rho)$ in the center direction and $b(g)$ measures the maximum contraction in $B(q, \rho)$ in the center direction.
We know that if $\pi_g(z)=\pi_g(y)$, then $d(z,y)<2\delta$.  So if $y\in \mathbb{T}^d-B(q,\rho+2\delta)$, then $z\notin B(q, \rho)$ and
$$|Dg_z\mathcal{F}^c|\geq a(g)\geq \lambda_c-\gamma.$$  Fix $\sigma>0$ such that
$$(\lambda_c-\gamma)^{1-m-\sigma}(b(f)-\gamma)^{2m+\sigma}>1.$$  Hence, for $\mu$-almost
every $x\in\mathbb T^d$, there exists some $K(x)>0$ such that, for all $z\in\pi_g^{-1}(x)$,
all $k\geq0$, and
$$\begin{array}{rlll}
|Dg^k_z\mathcal{F}^c| & \geq K(x)[a(g)^{1-m-\sigma} b(g)^{2m+\sigma}]^k\\
 & \geq K(x)[(\lambda_c-\gamma)^{1-m-\sigma}(b(f)-\gamma)^{2m+\sigma}]^k\\
  & \geq K(x)c^k
  \end{array}$$
with $c>1$. As $\pi_g^{-1}(x)$ must keep a bounded length it must be a unique point
for $\mu$-almost every $x$.  This shows that $\nu=\bar{\mu}(\textrm{mod }0)$ and $g$ is intrinsically ergodic. $\Box$

\bibliographystyle{plain}
\bibliography{maxentropy}

\end{document}